\documentclass{article}
\usepackage[margin=2.5cm]{geometry}
\usepackage{mathtools}
\usepackage{amssymb}
\usepackage{tikz}
\usepackage{algorithm}
\usetikzlibrary{decorations.markings, decorations.pathreplacing}
\usepackage[font=small, width=0.875\textwidth]{caption}
\numberwithin{equation}{section}
\makeatletter
\newenvironment{squareCases}{
  \matrix@check\squareCases\env@squareCases
}{
  \endarray\right.
}
\def\env@squareCases{
  \let\@ifnextchar\new@ifnextchar
  \left[
  \def\arraystretch{1.75}
  \array{@{}l@{}}
}
\begin{document}

\title{\huge Conjectures on union-closed families of sets}
\author{\Large Christopher Bouchard}
\date{}
\maketitle

\medskip

\abstract{\noindent A family of sets $\mathcal{A}$ is union-closed if it is finite and nonempty with member sets that are all finite and distinct (at least one of which is nonempty) and it satisfies the property $X, Y \in \mathcal{A} \implies X \cup Y \in \mathcal{A}$. Let $\binom{S}{k}$ be the set of all $k$-element subsets of a set $S$, and let $[n]=\{1,2,\cdots,n\}$ represent $\bigcup_{A \in \mathcal{A}}A$. Further, let $\mathcal{A}_B=\{A\in\mathcal{A} \ | \ A \cap B = B\}$ and $\mathcal{A}_{\underline{B}}=\{A\in\mathcal{A} \ | \ A \cap B = \emptyset\}$. We consider, for any union-closed family $\mathcal{A}$, the class of conjectures $\textrm{UC}_x \colon \ \exists B \in \binom{[n]}{n-x+1} \ | \ |\mathcal{A}_B| \geq |\mathcal{A}_{\underline{B}}|$, where $x \in [n]$. The extremal case $x=n$ is equivalent to the union-closed sets conjecture, also known as Frankl's conjecture, which states that there exists an element of $[n]$ that is in at least $\frac{|\mathcal{A}|}{2}$ member sets of $\mathcal{A}$. We prove $\textrm{UC}_x$ for $x \in [\lceil \frac{n}{3} \rceil + 1]$, and also investigate two strengthenings of the union-closed sets conjecture.}

\bigskip

\section*{1. Introduction}

A family of sets $\mathcal{A}$ is \textit{union-closed} if it is finite and nonempty with member sets that are all finite and distinct (at least one of which is nonempty) and it satisfies the property $X, Y \in \mathcal{A} \implies X \cup Y \in \mathcal{A}$. For a union-closed family $\mathcal{A}$, let $[n]=\{1,2,\cdots,n\}$ represent its \textit{base set} $\bigcup_{A \in \mathcal{A}}A$. The union-closed sets conjecture, also known as Frankl's conjecture, states the following:

\medskip

\noindent \textbf{Conjecture 1.1.} \textit{For any union-closed family} $\mathcal{A}$\textit{, there exists an element of} $[n]$ \textit{that is in at least} $\frac{|\mathcal{A}|}{2}$ \textit{member sets of} $\mathcal{A}$\textit{.}

\medskip

Conjecture 1.1 has been proved for $n \leq 12$ (see [11]) and $|\mathcal{A}| \leq 50$ (see [5] and [9] together with [11]). Equivalent conjectures include the lattice formulation of the union-closed sets conjecture, which has been verified for lower semimodular lattices (see [8]), and the graph formulation of the union-closed sets conjecture (see [1]). The \textit{frequency} in $\mathcal{A}$ of an element $y \in [n]$ is the number of member sets in $\mathcal{A}$ that contain $y$. A popular approach to Conjecture 1.1 is to improve the lower bound for the highest frequency of an element from $[n]$, proving a sequence of increasingly stronger results. In [3], entropic methods were used to prove that there exists an element from $[n]$ that is in at least $\varepsilon |\mathcal{A}|$ member sets of $\mathcal{A}$ for a constant $\varepsilon=\frac{1}{100}$, a result which in [10] was improved to $\varepsilon=\frac{3-\sqrt{5}}{2}$. Another well-studied approach aims to prove the full lower bound $\frac{|\mathcal{A}|}{2}$ for union-closed families that decrease in size $|\mathcal{A}|$ as a function of $n$. In [4], methods of Boolean analysis were applied to prove Conjecture 1.1 for $|\mathcal{A}| \geq 2^{n-1}$. For a survey of many results regarding Conjecture 1.1, see [2].

Let $\mathcal{A}_B=\{A\in\mathcal{A} \ | \ A \cap B = B\}$, $\mathcal{A}_{\underline{B}}=\{A\in\mathcal{A} \ | \ A \cap B = \emptyset\}$, and $\mathcal{A}_{B_1\underline{B_2}}=\mathcal{A}_{B_1} \cap \mathcal{A}_{\underline{B_2}}$. Also, let $\binom{S}{k}$ be the set of all $k$-element subsets of a set $S$. In the present work, we consider the following class of conjectures for any union-closed family $\mathcal{A}$:

\begin{equation}\label{UCx}
\underset{x \in [n]}{\operatorname{\textrm{UC}_x}} \colon \, \exists B \in \binom{[n]}{n-x+1} \textrm{ such that } |\mathcal{A}_B| \geq |\mathcal{A}_{\underline{B}}|.
\end{equation}

\noindent Conjecture 1.1 is itself the extremal case $x=n$ of $\textrm{UC}_x$, being equivalent to the statement $\exists B \in \binom{[n]}{1} \ | \ |\mathcal{A}_{B}| \geq |\mathcal{A}_{\underline{B}}|$. Therefore, another approach to Conjecture 1.1 is to prove $\textrm{UC}_x$ for increasing $x$. We prove $\textrm{UC}_x$ for $x \in [\lceil \frac{n}{3} \rceil + 1]$ and pose a question which, if answered in the affirmative, would prove $\textrm{UC}_x$ for $x \in [\lfloor \frac{n}{2} \rfloor + 1]$, where $\lceil \frac{n}{3} \rceil$ is the smallest integer greater than or equal to $\frac{n}{3}$ and $\lfloor \frac{n}{2} \rfloor$ is the largest integer less than or equal to $\frac{n}{2}$. We then proceed to investigate two strengthenings of Conjecture 1.1: one that imposes further constraint on the nested structure of union-closed families, and another that proposes an alternate definition of a finite power set.

\section*{2. Considering $\textrm{UC}_x$ for various values of $x \in [n]$}

We now prove, for any union-closed family $\mathcal{A}$, $\textrm{UC}_x$ of \eqref{UCx} whenever $x \in [\lceil \frac{n}{3} \rceil +1]$.

\medskip

\noindent \textbf{Theorem 2.1.} \textit{For any union-closed family} $\mathcal{A}$\textit{,} $\forall x \in [\lceil \frac{n}{3} \rceil +1] \ \textrm{UC}_x$\textit{.}

\medskip

\noindent \textit{Proof.} The case $n=1$ consists of the statements $\exists B \in \binom{[1]}{1} \ | \ |\mathcal{A}_B| \geq |\mathcal{A}_{\underline{B}}|$ and $\exists B \in \binom{[1]}{0} \ | \ |\mathcal{A}_B| \geq |\mathcal{A}_{\underline{B}}|$, for which $B=[1]$ and $B=\emptyset$, respectively. For $n>1$, we use induction on $x$.

\medskip

\noindent \textit{Base Case (}$x=1$\textit{):} $\exists B \in \binom{[n]}{n} \ | \ |\mathcal{A}_{B}| \geq |\mathcal{A}_{\underline{B}}|$.

\medskip

\noindent \textit{Proof.} $[n]=\bigcup_{A \in \mathcal{A}}A$ and $\mathcal{A}$ is union-closed. Therefore, $[n] \in \mathcal{A}$ and $|\mathcal{A}_{[n]}|=1$. If $\emptyset \in \mathcal{A}$ then $|\mathcal{A}_{\underline{[n]}}|=1$, and if $\emptyset \not \in \mathcal{A}$ then $|\mathcal{A}_{\underline{[n]}}|=0$. In either case, $B=[n] \in \binom{[n]}{n}$ with $|\mathcal{A}_{B}| \geq |\mathcal{A}_{\underline{B}}|$.

\medskip

\noindent \textit{Induction Step:} $\forall x \in [\lceil \frac{n}{3} \rceil ] (\exists B \in \binom{[n]}{n-x+1} \ | \ |\mathcal{A}_B| \geq |\mathcal{A}_{\underline{B}}| \implies \exists B' \in \binom{[n]}{n-x} \ | \ |\mathcal{A}_{B'}| \geq |\mathcal{A}_{\underline{B'}}|)$.

\medskip

\noindent \textit{Proof.} We recall the notation $\mathcal{A}_{B_1\underline{B_2}}=\mathcal{A}_{B_1} \cap \mathcal{A}_{\underline{B_2}}$. 

\medskip

\noindent If $\exists y' \in  B \ | \ \mathcal{A}_{\{y'\} \underline{B \setminus \{y'\}}} = \emptyset$, then $B' = B \setminus \{y'\}$ with $|\mathcal{A}_{B'}| \geq |\mathcal{A}_B| \geq |\mathcal{A}_{\underline{B}}| = |\mathcal{A}_{\underline{B'}}|$. 

\medskip

\noindent Else, $\forall y \in  B \ \mathcal{A}_{\{y\} \underline{B \setminus \{y\}}} \neq \emptyset$, in which case we have the following:

\medskip

\noindent $\exists C \in \binom{[n]}{n-x}$ in $\mathcal{A}$.

\medskip

\noindent \textit{Proof.} Consider the family $\mathcal{P}=\bigcup_{y \in B} \mathcal{A}_{\{y\} \underline{B \setminus \{y\}}}$ and the set $P=\bigcup_{p \in \mathcal{P}}p$. We have that $B \subseteq P \subseteq [n]$. By the construction of $P$, $\forall w \in P \setminus B \ \exists u \in B$ with a $z \in  \mathcal{A}_{\{u\} \underline{B \setminus \{u\}}}$ such that $w \in z$. Let $\mathcal{Z}$ be the family containing a (possibly repeated) $z$ for every $w$, so $|\mathcal{Z}| = |P \setminus B|$, and let $Z=\bigcup_{z \in \mathcal{Z}}z$. Then $Z \supseteq P \setminus B$ and $|Z| \leq 2|P \setminus B|$. By the union-closed property, $Z$ must be in $\mathcal{A}$ if $Z \neq \emptyset$. We observe that together $|B|=n-x+1$, $x \leq \lceil \frac{n}{3} \rceil$, $|P| \leq n$, and $|Z| \leq 2|P \setminus B|$ imply that $|Z| \leq 2 \lceil \frac{n}{3} \rceil - 2$. Similarly, $|B|=n-x+1$, $x \leq \lceil \frac{n}{3} \rceil$, and $|P| \geq |B|$ imply that $|P| \geq n-\lceil \frac{n}{3} \rceil + 1$. Therefore, $|Z| < |P|$, and $P \setminus Z \neq \emptyset$. Now, consider a $V \subsetneq P \setminus Z$ with $|V|=n-x-|Z|$. Such a $V$ must exist because of the following two facts: 

\medskip
 
\noindent 1.) $n-x-|Z| \geq 0$, as $x \leq \lceil \frac{n}{3} \rceil$ and $|Z|\leq 2|P \setminus B| \leq 2|[n]\setminus B| = 2(n-(n-x+1))=2x-2$, making $n-x-|Z|$ have a minimum of $n-3 \lceil \frac{n}{3} \rceil+2$, which itself has a minimum of $0$ when $n=3k+1$ for some whole number $k$. 

\medskip

\noindent 2.) $n-x-|Z| < |P \setminus Z|$, as $|P \setminus Z| = |P|-|Z|$ and $|P| \geq |B| = n-x+1$.

\medskip
 
\noindent Because $V \subsetneq P \setminus Z \subseteq B$ and $\forall y \in  B \ \mathcal{A}_{\{y\} \underline{B \setminus \{y\}}} \neq \emptyset$, we have that $\forall v \in V \ \exists V' \in \mathcal{A}_{\{v\} \underline{B \setminus \{v\}}}$. We collect a $V'$ for every $v$ into a family $\mathcal{V}$ and let $F=\bigcup_{V' \in \mathcal{V}}V'$. By the union-closed property, $F$ must be in $\mathcal{A}$ if $F \neq \emptyset$. Also, $V \subseteq F$ and $F \setminus V \subseteq P \setminus B \subseteq Z$, making $Z \cup F=Z \cup V$. And because $Z \cap V = \emptyset$, we have $|Z \cup V|=|Z|+|V|=|Z|+(n-x-|Z|)=n-x$. Thus, $Z \cup F \in \binom{[n]}{n-x}$, and, noting that $Z$ and $F$ cannot both be empty, we have that $Z \cup F$ is in $\mathcal{A}$ by the union-closed property. 

\medskip 

\noindent Therefore, $C = Z \cup F$, and $C \in \binom{[n]}{n-x}$ is in $\mathcal{A}$.

\medskip

\noindent It follows that $|\mathcal{A}_C| \geq |\mathcal{A}_{\underline{C}}|$, as any $X$ in $\mathcal{A}_{\underline{C}}$ can be uniquely matched with $C \cup X$ in $\mathcal{A}_C$. $C \cup X$ must be in $\mathcal{A}_C$ by the existence of $C$ in $\mathcal{A}$ and the union-closed property, and every $C \cup X$ must be unique as $\forall X \in \mathcal{A}_{\underline{C}} \ C \cap X = \emptyset$.

\medskip 

\noindent Thus, $B'=C$ is in $\binom{[n]}{n-x}$ with $|\mathcal{A}_{B'}| \geq |\mathcal{A}_{\underline{B'}}|$.
 
\medskip 

\noindent The induction step is proved, as we have proved the existence of a $B'$ in $\binom{[n]}{n-x}$ such that $|\mathcal{A}_{B'}| \geq |\mathcal{A}_{\underline{B'}}|$ when $\exists y' \in  B \ | \ \mathcal{A}_{\{y'\} \underline{B \setminus \{y'\}}} = \emptyset$, and when $\forall y \in  B \ \mathcal{A}_{\{y\} \underline{B \setminus \{y\}}} \neq \emptyset$.

\medskip

\noindent As the induction step is proved, the proof of Theorem 2.1 is complete.

\medskip 

The proof method reached its limit as it achieved a minimum of $n-x-|Z|=0$. We now consider a question which, if answered in the affirmative, would extend the result to $x \in [\lfloor \frac{n}{2} \rfloor +1]$.

\medskip 

\noindent \textbf{Question 2.2.} If $\mathcal{F}$ is a (not necessarily union-closed) finite family of (not necessarily distinct) finite sets such that $|\mathcal{F}| > |\bigcup_{F \in \mathcal{F}}F|+1$, must there exist $\mathcal{F}'\subsetneq \mathcal{F}$ with $|\mathcal{F}| - |\mathcal{F}'|=|\bigcup_{F \in \mathcal{F}'}F|+1$?

\medskip 

\noindent \textbf{Proposition 2.3.} \textit{An affirmative answer to Question 2.2 would prove, for any union-closed family} $\mathcal{A}$\textit{,} $\textrm{UC}_x$ \textit{for all} $x \in [\lfloor \frac{n}{2} \rfloor +1]$\textit{.}

\medskip 

\noindent \textit{Proof.} The proof for $x \in [\lfloor \frac{n}{2} \rfloor +1]$ would start the same as the proof for $x \in [\lceil \frac{n}{3} \rceil +1]$ started (except for the case $n=1$ now consisting only of the statement $\exists B \in \binom{[1]}{1} \ | \ |\mathcal{A}_B| \geq |\mathcal{A}_{\underline{B}}|$ (for which $B=[1]$), and the induction step for $n>1$ being extended from $x \in [\lceil \frac{n}{3}\rceil]$ to $x \in [\lfloor \frac{n}{2} \rfloor]$). The proofs differ in the induction step in how they show existence of $C$. We resume the proof at that point:
 
\medskip

\noindent $\exists C \in \binom{[n]}{n-x}$ in $\mathcal{A}$.

\medskip 

\noindent \textit{Proof.} Recall that at this point in the proof we had assumed that $\forall y \in  B \ \mathcal{A}_{{\{y\}}{\underline{B \setminus \{y\}}}} \neq \emptyset$. Therefore, there must be a family $\mathcal{S} \subsetneq \mathcal{A}$ with $|\mathcal{S}|=|B|$ such that $\forall b \in B \ \exists S \in \mathcal{S} \cap \mathcal{A}_{\{b\} \underline{B \setminus \{b\}}}$. 

\medskip

\noindent In this case, $\mathcal{F}$ from Question 2.2 corresponds to $\mathcal{T}=\{S \cap ([n] \setminus B) \ | \ S \in \mathcal{S}\}$. $\mathcal{T}$ satisfies the definition of $\mathcal{F}$ because $|\mathcal{T}| > |\bigcup_{T \in \mathcal{T}}T|+1$, as $|\mathcal{T}|=|B|=n-x+1 \geq n-\lfloor \frac{n}{2} \rfloor +1$ and $|\bigcup_{T \in \mathcal{T}}T| \leq |[n] \setminus B|=x-1 \leq \lfloor \frac{n}{2} \rfloor - 1$. Further, $\mathcal{T}$ need not be union-closed, and could have repeated member sets if some of the member sets of $\mathcal{S}$ shared exactly the same elements from $[n] \setminus B$.

\medskip

\noindent Now applying Question 2.2 to $\mathcal{T}$, there exists $\mathcal{T}' \subsetneq \mathcal{T}$ with $|\mathcal{T}| - |\mathcal{T}'|=|\bigcup_{T \in \mathcal{T}'}T|+1$. It follows that $|\bigcup_{T \in \mathcal{T}'}T|+|\mathcal{T}'|=|\mathcal{T}|-1=|B|-1=(n-x+1)-1=n-x$. Recall that every member set $T$ of $\mathcal{T}$ could be written as $T= S\setminus \{b\}$ for some unique $S \in \mathcal{S}$ and $b \in B$. Let $\mathcal{S}' \subsetneq \mathcal{S}$ be the family consisting of the corresponding $S$ for every $T \in \mathcal{T}'$. Then $|\bigcup_{S \in \mathcal{S}'}S|=|\bigcup_{T \in \mathcal{T}'}T|+|\mathcal{T}'|=n-x$, as $\bigcup_{S \in \mathcal{S}'}S=(\bigcup_{T \in \mathcal{T}'}T) \cup B^*$ and $(\bigcup_{T \in \mathcal{T}'}T) \cap B^*=\emptyset$, where $B^* \in \binom{B}{|\mathcal{T}'|}$. Finally, $\bigcup_{S \in \mathcal{S}'}S \in \mathcal{A}$ by $\mathcal{S}'\subsetneq \mathcal{A}$ and the union-closed property. Therefore, $C=\bigcup_{S \in \mathcal{S}'}S$. 

\medskip 

\noindent After $C$ is shown to exist, the proof for $x \in [\lfloor \frac{n}{2} \rfloor +1]$ is the same as that for $x \in [\lceil \frac{n}{3} \rceil +1]$. This concludes the proof of Proposition 2.3.

\medskip

We note that this particular method that finds $C \in \binom{[n]}{n-x}$ in $\mathcal{A}$ breaks down for $x > \lfloor \frac{n}{2} \rfloor$ in the induction step. This is because the method is not able to guarantee that $[n] \setminus B$ is not a proper subset of $C$, and so a fundamental problem arises when $|[n] \setminus B| \geq |C|$, which occurs in the induction step when $x \geq \frac{n+1}{2}$. 

\medskip

The following theorem also deals with $\textrm{UC}_x$ of \eqref{UCx}.

\medskip 
 
\noindent \textbf{Theorem 2.4.} \textit{For any union-closed family} $\mathcal{A}$ \textit{with} $n>1$\textit{,} $\textrm{UC}_{n-1} \implies \textrm{UC}_n$\textit{.} 

\medskip 

\noindent \textit{Proof.} Without loss of generality, let $\{1,2\}$ be the doubleton from $\textrm{UC}_{n-1}$, so $|\mathcal{A}_{\{1,2\}}| \geq |\mathcal{A}_{\underline{\{1,2\}}}|$. Using the notation $\mathcal{A}_{\{1\}\underline{\{2\}}}=\mathcal{A}_{\{1\}} \cap \mathcal{A}_{\underline{\{2\}}}$ and $\mathcal{A}_{\{2\}\underline{\{1\}}}=\mathcal{A}_{\{2\}} \cap \mathcal{A}_{\underline{\{1\}}}$, we observe that if $|\mathcal{A}_{\{1\}\underline{\{2\}}}| \geq |\mathcal{A}_{\{2\}\underline{\{1\}}}|$, then the singleton from $\textrm{UC}_n$ is equal to $\{1\}$, as $(\mathcal{A}_{\{1,2\}} \cup \mathcal{A}_{\{1\}\underline{\{2\}}}=\mathcal{A}_{\{1\}}) \land (\mathcal{A}_{\{1,2\}} \cap \mathcal{A}_{\{1\}\underline{\{2\}}}=\emptyset) \implies |\mathcal{A}_{\{1\}}|=|\mathcal{A}_{\{1, 2\}}|+|\mathcal{A}_{\{1\}\underline{\{2\}}}|$ and $(\mathcal{A}_{\underline{\{1,2\}}} \cup \mathcal{A}_{\{2\}\underline{\{1\}}}=\mathcal{A}_{\underline{\{1\}}}) \land (\mathcal{A}_{\underline{\{1,2\}}} \cap \mathcal{A}_{\{2\}\underline{\{1\}}} = \emptyset) \implies |\mathcal{A}_{\underline{\{1\}}}|=|\mathcal{A}_{\underline{{\{1, 2\}}}}|+|\mathcal{A}_{\{2\}\underline{\{1\}}}|$. Else, $|\mathcal{A}_{\{2\}\underline{\{1\}}}| > |\mathcal{A}_{\{1\}\underline{\{2\}}}|$, making the singleton from $\textrm{UC}_n$ equal to $\{2\}$, as $(\mathcal{A}_{\{1,2\}} \cup \mathcal{A}_{\{2\}\underline{\{1\}}}=\mathcal{A}_{\{2\}}) \land (\mathcal{A}_{\{1,2\}} \cap \mathcal{A}_{\{2\}\underline{\{1\}}}=\emptyset) \implies |\mathcal{A}_{\{2\}}|=|\mathcal{A}_{\{1, 2\}}|+|\mathcal{A}_{\{2\}\underline{\{1\}}}|$ and $(\mathcal{A}_{\underline{\{1,2\}}} \cup \mathcal{A}_{\{1\}\underline{\{2\}}}=\mathcal{A}_{\underline{\{2\}}}) \land (\mathcal{A}_{\underline{\{1,2\}}} \cap \mathcal{A}_{\{1\}\underline{\{2\}}} = \emptyset) \implies |\mathcal{A}_{\underline{\{2\}}}|=|\mathcal{A}_{\underline{{\{1, 2\}}}}|+|\mathcal{A}_{\{1\}\underline{\{2\}}}|$. This concludes the proof of Theorem 2.4. 

\medskip

$\textrm{UC}_{n-1} \implies \textrm{UC}_n$ in fact holds for any finite family of finite sets with $n>1$, whether the family be union-closed or not. Further, by considering the contrapositive, we have a necessary condition for any counterexample to Conjecture 1.1, namely that a counterexample $\tilde{\mathcal{A}}$ on a base set $[\tilde{n}]$ must have $|\tilde{\mathcal{A}}_{\{x,y\}}| < |\tilde{\mathcal{A}}_{\underline{\{x,y\}}}|$ for every $\{x,y\}$ in $\binom{[\tilde{n}]}{2}$. Another question is the following:

\medskip 

\noindent \textbf{Question 2.5.} For a union-closed family $\mathcal{A}$ with $n>1$, does $\textrm{UC}_n$ imply $\textrm{UC}_x$ for all $x \in [n-1]$?

\medskip 

We have already verified by Theorem 2.4 that for a union-closed family $\mathcal{A}$ with $n>1$, $\textrm{UC}_{n-1} \implies \textrm{UC}_n$. For answering Question 2.5, a starting point could be to prove the converse, i.e. that $\textrm{UC}_n \implies \textrm{UC}_{n-1}$. At the core of understanding Conjecture 1.1 is understanding the union-closed property in general. The following could help in this development:

\medskip 

\noindent 1.) Answering Question 2.2 in the affirmative, thus proving $\textrm{UC}_x$ for $x \in [\lfloor \frac{n}{2} \rfloor +1]$.

\medskip 

\noindent 2.) Proving certain implications with regard to $\textrm{UC}_x$, especially, for $n>4$, finding $x \in [n-1] \setminus [\lceil \frac{n}{3} \rceil+1]$ such that $\textrm{UC}_{x+1} \implies \textrm{UC}_x$.  

\section*{3. A strengthening of Conjecture 1.1}

Before presenting the strengthening of Conjecture 1.1, we discuss a theorem of Reimer in the context of an important characteristic of union-closed families. The binary logarithm and natural logarithm of a positive real number $x$ are denoted by $\log(x)$ and $\ln(x)$, respectively.

\medskip

\noindent \textbf{Theorem 3.1} (Reimer [7])\textbf{.} \textit{For any union-closed family} $\mathcal{A}$\textit{,} $\sum_{k=1}^{n} |\mathcal{A}_{\{k\}}| \ \geq \ \frac{|\mathcal{A}|}{2}\log(|\mathcal{A}|)$\textit{.}

\medskip 

Theorem 3.1 was proved conclusively by Reimer without need of assuming Conjecture 1.1. However, we show that Theorem 3.1 follows directly if we assume Conjecture 1.1, in order to highlight the nested structure of union-closed families.

\medskip 

\noindent \textbf{Proposition 3.2.} \textit{Conjecture 1.1} $\implies$ \textit{Theorem 3.1.}

\medskip 

\noindent \textit{Proof.} We show Theorem 3.1, assuming Conjecture 1.1. We use induction on size of the base set. 

\medskip 

\noindent \textit{Base Case (}$n=1$\textit{):} The only two union-closed families on the base set $[1]$ are $\{\emptyset, \{1\}\}$ and $\{\{1\}\}$. In the first case, $\sum_{k=1}^{n} |\mathcal{A}_{\{k\}}|=1 \geq \frac{|\mathcal{A}|}{2}\log(|\mathcal{A}|)=\log(2)=1$. In the second case, $\sum_{k=1}^{n} |\mathcal{A}_{\{k\}}|=1 \geq \frac{|\mathcal{A}|}{2}\log(|\mathcal{A}|)=\frac{1}{2}\log(1)=0$.

\medskip 

\noindent \textit{Induction Step:} Without loss of generality, we consider a union-closed family $\mathcal{A}$ with $n > 1$ such that $|\mathcal{A}_{\{n\}}|=\max_{x \in [n]} \{|\mathcal{A}_{\{x\}}|\}$. We show that:
\noindent $\Bigr($For every union-closed family $\mathcal{A}^*$ on a base set $N^* \subseteq [n-1]$, $\sum_{k \in N^*} |{\mathcal{A}^*}_{\{k\}}| \geq \frac{|\mathcal{A}^*|}{2}\log(|\mathcal{A}^*|)\Bigr) \implies \sum_{k=1}^n |\mathcal{A}_{\{k\}}| \geq \frac{|\mathcal{A}|}{2}\log(|\mathcal{A}|)$.

\medskip 

\noindent \textit{Proof.} We have $c =\frac{|\mathcal{A}_{\{n\}}|}{|\mathcal{A}|} \in [\frac{1}{2}, 1]$, as we are assuming Conjecture 1.1. Let ${\mathcal{A}_{\{n\}}}'=\{A \setminus \{n\} \ | \ A \in \mathcal{A}_{\{n\}}\}$. It follows that $\bigcup_{A \in {{\mathcal{A}_{\{n\}}}}'}A=[n-1]$, and $\bigcup_{A \in \mathcal{A}_{\underline{\{n\}}}} A = N \subseteq [n-1]$. We consider three cases:

\medskip

\noindent 1.) $\mathcal{A}_{\underline{\{n\}}}=\emptyset$ (i.e. $c=1$): By the hypothesis of the induction step, we have $\sum_{k=1}^{n-1} |{{(\mathcal{A}_{\{n\}}}')}_{\{k\}}| \geq \frac{|{\mathcal{A}_{\{n\}}}'|}{2} \log(|{\mathcal{A}_{\{n\}}}'|)$. In this case, $|{\mathcal{A}_{\{n\}}}'|=|\mathcal{A}|$, and $\sum_{k=1}^{n-1} |{({\mathcal{A}_{\{n\}}}')}_{\{k\}}|=(\sum_{k=1}^{n} |\mathcal{A}_{\{k\}}|)-|\mathcal{A}|$. Therefore, $\sum_{k=1}^{n} |\mathcal{A}_{\{k\}}| \geq \frac{|\mathcal{A}|}{2}\log(|\mathcal{A}|) + |\mathcal{A}| > \frac{|\mathcal{A}|}{2}\log(|\mathcal{A}|)$.

\medskip
\smallskip

\noindent 2.) $\mathcal{A}_{\underline{\{n\}}}=\{\emptyset\}$: By the hypothesis of the induction step, we again have $\sum_{k=1}^{n-1} |{{(\mathcal{A}_{\{n\}}}')}_{\{k\}}| \geq \frac{|{\mathcal{A}_{\{n\}}}'|}{2} \log(|{\mathcal{A}_{\{n\}}}'|)$. In this case, $|{\mathcal{A}_{\{n\}}}'|=|\mathcal{A}|-1$, and $\sum_{k=1}^{n-1} |{({\mathcal{A}_{\{n\}}}')}_{\{k\}}|=(\sum_{k=1}^{n} |\mathcal{A}_{\{k\}}|)-(|\mathcal{A}|-1)$. Therefore, $\sum_{k=1}^{n} |\mathcal{A}_{\{k\}}| \geq \frac{|\mathcal{A}|-1}{2} \log (|\mathcal{A}|-1) + |\mathcal{A}| -1$. If $\frac{|\mathcal{A}|}{2}\log(|\mathcal{A}|)$ is a lower bound for $\frac{|\mathcal{A}|-1}{2} \log (|\mathcal{A}|-1) + |\mathcal{A}| -1$, then it must also be a lower bound for $\sum_{k=1}^{n} |\mathcal{A}_{\{k\}}|$. Thus, in order to prove that $\sum_{k=1}^n |\mathcal{A}_{\{k\}}| \geq \frac{|\mathcal{A}|}{2}\log(|\mathcal{A}|)$, it suffices to show that $\frac{|\mathcal{A}|-1}{2} \log (|\mathcal{A}|-1) + |\mathcal{A}| -1 \geq \frac{|\mathcal{A}|}{2}\log(|\mathcal{A}|)$. This is equivalent to showing that $\left(\frac{2|\mathcal{A}|-2}{|\mathcal{A}|}\right)^{|\mathcal{A}|-1}\left(\frac{2^{|\mathcal{A}|-1}}{|\mathcal{A}|}\right) \geq 1$, which is true when $|\mathcal{A}| \geq 2$. Now, $\mathcal{A}_{\underline{\{n\}}}=\{\emptyset\} \implies \emptyset \in \mathcal{A}$, and because $\mathcal{A}$ is union-closed, $\mathcal{A}$ must also contain its base set $[n]$. Thus, $|\mathcal{A}| \geq 2$ and $\sum_{k=1}^n |\mathcal{A}_{\{k\}}| \geq \frac{|\mathcal{A}|}{2}\log(|\mathcal{A}|)$.

\medskip
\smallskip

\noindent 3.) $\mathcal{A}_{\underline{\{n\}}} \neq \emptyset \ \land \ \mathcal{A}_{\underline{\{n\}}} \neq \{\emptyset\}$: In this case, $\mathcal{A}_{\underline{\{n\}}}$ is union-closed in addition to ${\mathcal{A}_{\{n\}}}'$ being union-closed, and we establish a lower bound for $\sum_{k=1}^{n} |\mathcal{A}_{\{k\}}|$ by solving the following problem:

\[\min\limits_{c \in [\frac{1}{2},1)}\{ I(c) \}\textrm{, where }I(c)=c|\mathcal{A}| + \frac{c|\mathcal{A}|}{2}\log(c|\mathcal{A}|)+\frac{(1-c)|\mathcal{A}|}{2}\log((1-c)|\mathcal{A}|)\textrm{.}\]

\noindent In $I(c)$, the first term $c|\mathcal{A}|$ comes from the $c|\mathcal{A}|=|\mathcal{A}_{\{n\}}|$ number of times that element $n$ was removed from $\mathcal{A}_{\{n\}}$ to make ${\mathcal{A}_{\{n\}}}'$. By the hypothesis of the induction step, the second term $\frac{c|\mathcal{A}|}{2}\log(c|\mathcal{A}|)$ is a lower bound of $\sum_{k=1}^{n-1}|{({\mathcal{A}_{\{n\}}}')}_{\{k\}}|$, as $c|\mathcal{A}|=|\mathcal{A}_{\{n\}}|=|{\mathcal{A}_{\{n\}}}'|$, and the final term $\frac{(1-c)|\mathcal{A}|}{2}\log((1-c)|\mathcal{A}|)$ is a lower bound of $\sum_{k \in N}|{(\mathcal{A}_{{\underline{\{n\}}}})}_{\{k\}}|$, as $(1-c)|\mathcal{A}|=|\mathcal{A}_{{\underline{\{n\}}}}|$. For $c \in [\frac{1}{2}, 1)$, $\frac{dI}{dc}=|\mathcal{A}|+\frac{|\mathcal{A}|}{2}(\log(c)-\log(1-c))>0$, as $\frac{dI}{dc} \Bigr|_{c=\frac{1}{2}}=|\mathcal{A}|>0$ and when $c \in [\frac{1}{2},1)$, $\frac{d^2I}{dc^2}=\frac{|\mathcal{A}|}{2\ln(2)}(\frac{1}{c(1-c)})>0$. Thus, the global minimum of $I(c)$ for $c \in [\frac{1}{2},1)$ is achieved uniquely at $c=\frac{1}{2}$. Plugging this minimizer into $I(c)$ gives $I(\frac{1}{2})=(\frac{1}{2})|\mathcal{A}| + \frac{(\frac{1}{2})|\mathcal{A}|}{2}\log((\frac{1}{2})|\mathcal{A}|)+\frac{(1-\frac{1}{2})|\mathcal{A}|}{2}\log((1-\frac{1}{2})|\mathcal{A}|)=\frac{|\mathcal{A}|}{2}\log(|\mathcal{A}|)$, which is exactly the lower bound that we intended to prove for $\sum_{k=1}^n |\mathcal{A}_{\{k\}}|$. This proves the final case of the induction step, concluding the proof of Proposition 3.2.

\medskip

In general, induction can be useful for proving properties of union-closed families, especially when coupled with additional assumptions, because union-closed families are rich in subfamilies that are also union-closed. The nested structure of union-closed families motivates the following strengthening of Conjecture 1.1.

\medskip

\noindent \textbf{Conjecture 3.3.} \textit{For any union-closed family} $\mathcal{A}$\textit{, if} $y \in [n]$ \textit{with} $|\mathcal{A}_{\{y\}}| = \max_{x \in [n]}\{|\mathcal{A}_{\{x\}}|\}$ \textit{and there exists} $z \in [n]$ \textit{such that} $\mathcal{A} = \mathcal{A}_{\{y\}} \cup \mathcal{A}_{\{z\}}$\textit{, then} $|\mathcal{A}_{\{y\}}| \geq 2|\mathcal{A}_{\underline{\{y\}}}|$\textit{.}

\medskip

Before proving that Conjecture 3.3 implies Conjecture 1.1, we present a framework for considering union-closed families. For any union-closed family $\mathcal{A}$, we define the following:

\medskip
\smallskip

\begin{center}

\def\arraystretch{1.5}

$\begin{array}{|l}
\mathcal{R}^{(i)} \colon i \in \mathbb{Z}_{\geq 0}\textrm{,}
\\
\mathcal{P}=\{\mathcal{P}^{(i)} \ | \ i \in [\min(\{ i \ | \ \mathcal{R}^{(i)} = \emptyset \})]\}\textrm{, and}
\\
P \colon [\min(\{ i \ | \ \mathcal{R}^{(i)} = \emptyset \})] \to [n]\textrm{,}
\\
\textrm{such that }\mathcal{R}^{(0)}=\mathcal{A}\textrm{, and:} 
\\
\textrm{For }i \in \mathbb{Z}_{>0}\textrm{:}
\\
\begin{squareCases}
\textrm{If }\mathcal{R}^{(i-1)} \neq \emptyset\textrm{, then:} 
\\
\begin{squareCases}
P(i)=\min\Bigr(\Bigr\{y \ \Bigr| \ |(\mathcal{R}^{(i-1)})_{\{y\}}|=\max_{x \in \bigcup_{R \in \mathcal{R}^{(i-1)}}R}\{|(\mathcal{R}^{(i-1)})_{\{x\}}|\}\Bigr\}\Bigr)\textrm{.}
\\
\textrm{If }(\mathcal{R}^{(i-1)})_{\underline{\{P(i)\}}} \neq \{\emptyset\}\textrm{, then }\mathcal{P}^{(i)}=(\mathcal{R}^{(i-1)})_{\{P(i)\}}\textrm{ and }\mathcal{R}^{(i)}=(\mathcal{R}^{(i-1)})_{\underline{\{P(i)\}}}\textrm{.}
\\
\textrm{If }(\mathcal{R}^{(i-1)})_{\underline{\{P(i)\}}} = \{\emptyset\}\textrm{, then }\mathcal{P}^{(i)}=(\mathcal{R}^{(i-1)})_{\{P(i)\}} \cup \{\emptyset\}\textrm{ and }\mathcal{R}^{(i)}=\emptyset\textrm{.}
\end{squareCases}
\\
\textrm{If }\mathcal{R}^{(i-1)} = \emptyset\textrm{, then }\mathcal{R}^{(i)}=\emptyset\textrm{.} 
\end{squareCases}
\end{array}$

\end{center}

\medskip
\smallskip

We note that $\mathcal{P}$ is a partition of $\mathcal{A}$ with $|\mathcal{P}|=\min(\{ i \ | \ \mathcal{R}^{(i)} = \emptyset \})$.

\medskip

\noindent \textbf{Theorem 3.4.} \textit{If} $S$ \textit{is a nonempty subset of} $[|\mathcal{P}|]$\textit{, then} $\bigcup_{k \in S} \mathcal{P}^{(k)}$ \textit{is union-closed.}

\medskip

\noindent \textit{Proof.} $\bigcup_{k \in S} \mathcal{P}^{(k)}$ is a subfamily of $\mathcal{A}$, so all of its member sets are finite and distinct. Also, $\bigcup_{k \in S} \mathcal{P}^{(k)}$ has at least one nonempty member set. To prove the union-closed property for $\bigcup_{k \in S} \mathcal{P}^{(k)}$, we observe that $(X,Y \in \bigcup_{k \in S} \mathcal{P}^{(k)}) \implies (X \in \mathcal{P}^{(i)} \land Y \in \mathcal{P}^{(l)})$, where $1 \leq i \leq l \leq |\mathcal{P}|$ without loss of generality. $X \cup Y \in \mathcal{A}$ because $X,Y \in \mathcal{A}$ and $\mathcal{A}$ is union-closed. If $X = \emptyset$, then $X \cup Y = Y \in \mathcal{P}^{(l)} \subseteq \bigcup_{k \in S} \mathcal{P}^{(k)}$. Else, $P(i)$ is in $X$, implying that $P(i)$ is also in $X \cup Y$. If $i=1$, then $(X \cup Y \in \mathcal{A}) \land (P(i) \in X \cup Y) \implies X \cup Y \in \mathcal{P}^{(i)} \subseteq \bigcup_{k \in S} \mathcal{P}^{(k)}$. If $i>1$, then $\forall j \in [i-1] \ ((P(j) \not \in X \land P(j) \not \in Y) \implies P(j) \not \in X \cup Y)$, and it follows that $((X \cup Y \in \mathcal{A}) \land (P(i) \in X \cup Y) \land (\forall j \in [i-1] \ P(j) \not \in X \cup Y)) \implies X \cup Y \in \mathcal{P}^{(i)} \subseteq \bigcup_{k \in S} \mathcal{P}^{(k)}$. Therefore, $\bigcup_{k \in S} \mathcal{P}^{(k)}$ satisfies the union-closed property, and the proof of Theorem 3.4 is complete.

\medskip

\noindent \textbf{Corollary 3.5.} \textit{Conjecture 1.1 implies that} $|\mathcal{P}^{(j)}| \geq \sum_{i \in S }|\mathcal{P}^{(i)}|$ \textit{for any nonempty subset} $S$ \textit{of} $[|\mathcal{P}|] \setminus [k]$\textit{, whenever }$1 \leq j \leq k < |\mathcal{P}|$\textit{.}

\medskip

\noindent \textbf{Theorem 3.6.} \textit{Conjecture 3.3} $\implies$ \textit{Conjecture 1.1.}

\medskip

\noindent \textit{Proof.} We consider a union-closed family $\mathcal{A}$ with $|\mathcal{P}|>1$, as Conjecture 1.1 holds trivially when $|\mathcal{P}|=1$. Assume that $\emptyset \not \in \mathcal{A}$. Theorem 3.6 will be proved when we show that Conjecture 1.1 holds for both $\mathcal{A}$ and $\mathcal{A} \cup \{\emptyset\}$. For all $i$ in $[|\mathcal{P}|-1]$, $\mathcal{P}^{(i)} \cup \mathcal{P}^{(i+1)}$ is a union-closed family that, according to the framework under consideration, has its own partition $\hat{\mathcal{P}}=\{\hat{\mathcal{P}}^{(1)}, \hat{\mathcal{P}}^{(2)}\}$ and function $\hat{P} \colon \{1,2\} \to \{\hat{P}(1), \hat{P}(2)\}$ such that $\hat{\mathcal{P}}^{(1)}=\mathcal{P}^{(i)}$, $\hat{\mathcal{P}}^{(2)}=\mathcal{P}^{(i+1)}$, $\hat{P}(1)=P(i)$, and $\hat{P}(2)=P(i+1)$. Conjecture 3.3 then applies to $\mathcal{P}^{(i)} \cup \mathcal{P}^{(i+1)}$ with $y$ and $z$ from Conjecture 3.3 respectively equal to $\hat{P}(1)$ and $\hat{P}(2)$. Thus, $|\mathcal{P}^{(i)}|\geq 2|\mathcal{P}^{(i+1)}|$ for all $i \in [|\mathcal{P}|-1]$, implying the following upper bound $\mu(\mathcal{P})$ of $\sum_{i=2}^{|\mathcal{P}|} |\mathcal{P}^{(i)}|$:

\[\mu(\mathcal{P})=\overbrace{\left \lfloor{\frac{1}{2}|\mathcal{P}^{(1)}|}\right \rfloor + \left \lfloor{\frac{1}{2}\left \lfloor{\frac{1}{2}|\mathcal{P}^{(1)}|}\right \rfloor }\right \rfloor  + \left \lfloor{\frac{1}{2}\left \lfloor{\frac{1}{2}\left \lfloor{\frac{1}{2}|\mathcal{P}^{(1)}|}\right \rfloor}\right \rfloor }\right \rfloor+\cdots + \left \lfloor{\frac{1}{2}\left \lfloor{\frac{1}{2}\left \lfloor{\frac{1}{2}\cdots \left \lfloor{\frac{1}{2}|\mathcal{P}^{(1)}|}\right \rfloor \cdots}\right \rfloor }\right \rfloor }\right \rfloor}^{|\mathcal{P}|-1\text{ terms}}.\]

\smallskip

\[\sum_{i=2}^{|\mathcal{P}|} |\mathcal{P}^{(i)}| \leq \mu(\mathcal{P}) \leq |\mathcal{P}^{(1)}| \sum_{i=1}^{|\mathcal{P}|-1} \frac{1}{2^i} < |\mathcal{P}^{(1)}| \sum_{i=1}^{\infty} \frac{1}{2^i}=|\mathcal{P}^{(1)}|.\]

\medskip

\noindent Therefore, $|\mathcal{P}^{(1)}| > \sum_{i=2}^{|\mathcal{P}|} |\mathcal{P}^{(i)}|$, and we have that $P(1) \in [n]$ with $|\mathcal{A}_{\{P(1)\}}| = |\mathcal{P}^{(1)}| > \sum_{i=2}^{|\mathcal{P}|} |\mathcal{P}^{(i)}| = |\mathcal{A}_{\underline{\{P(1)\}}}|$. Conjecture 1.1 then holds for $\mathcal{A}$, as $|\mathcal{A}_{\{P(1)\}}| > \frac{|\mathcal{A}|}{2}$. 

\medskip

\noindent Now consider $\mathcal{A} \cup \{\emptyset\}$. We have that $(\mathcal{A} \cup \{\emptyset\})_{\{P(1)\}} = \mathcal{A}_{\{P(1)\}}$ and $(\mathcal{A} \cup \{\emptyset\})_{\underline{\{P(1)\}}} = \mathcal{A}_{\underline{\{P(1)\}}} \cup \{\emptyset\}$. It follows that $|(\mathcal{A \cup \{\emptyset\}})_{\{P(1)\}}| > |(\mathcal{A \cup \{\emptyset\}})_{\underline{\{P(1)\}}}|-1$, and so $|(\mathcal{A \cup \{\emptyset\}})_{\{P(1)\}}| \geq |(\mathcal{A \cup \{\emptyset\}})_{\underline{\{P(1)\}}}|$ because $|(\mathcal{A \cup \{\emptyset\}})_{\{P(1)\}}|$ and $|(\mathcal{A \cup \{\emptyset\}})_{\underline{\{P(1)\}}}|$ are both integers. Thus, $|(\mathcal{A} \cup \{\emptyset\})_{\{P(1)\}}| \geq \frac{|\mathcal{A} \cup \{\emptyset\}|}{2}$, and Conjecture 1.1 also holds for $\mathcal{A} \cup \{\emptyset\}$. This completes the proof of Theorem 3.6.

\medskip

For any union-closed family $\mathcal{A}$ with $|\mathcal{P}|>1$, we have that $|\mathcal{A}_{\{P(1)\}\underline{\{P(2)\}}}| \geq |\mathcal{A}_{\{P(2)\}\underline{\{P(1)\}}}|$, where $\mathcal{A}_{\{P(1)\}\underline{\{P(2)\}}}=\mathcal{A}_{\{P(1)\}} \cap \mathcal{A}_{\underline{\{P(2)\}}}$ and $\mathcal{A}_{\{P(2)\}\underline{\{P(1)\}}}=\mathcal{A}_{\{P(2)\}} \cap \mathcal{A}_{\underline{\{P(1)\}}}$. We observe that Conjecture 1.1 holds if $|\mathcal{P}|=2$. A proof technique for Conjecture 1.1 would show, for any union-closed family $\mathcal{A}$ with $|\mathcal{P}|>2$, that $|\mathcal{A}_{\{P(1),P(2)\}}| \geq |\mathcal{A}_{\underline{\{P(1), P(2)\}}}|$ (see Theorem 2.4). On the other hand, a proof technique for Conjecture 3.3 would show, for any union-closed family $\mathcal{A}$ with $\emptyset \not \in \mathcal{A}$ and $|\mathcal{P}|=2$, that, without loss of generality, $|\mathcal{A}_{\{P(1),P(2)\}}| \geq |\mathcal{A}_{\{P(2)\}\underline{\{P(1)\}}}|$. Figures 3.1 and 3.2 illustrate the respective techniques.

\medskip
\smallskip

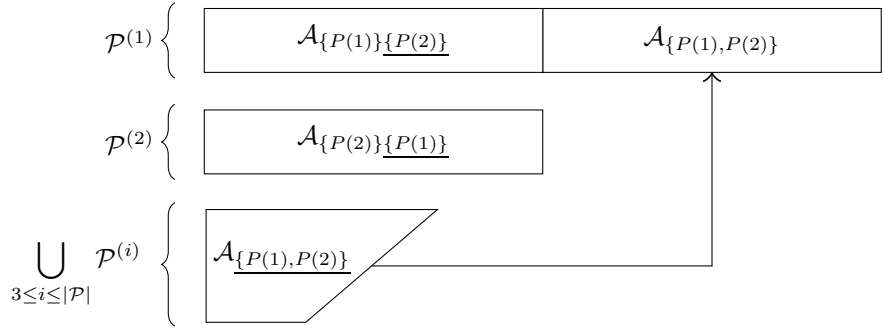
\begin{figure}[H]

\caption*{\textbf{Figure 3.1:} A technique for Conjecture 1.1 attempts to find, for any union-closed family $\mathcal{A}$ with $|\mathcal{P}|>2$, a unique $y \in \mathcal{A}_{\{P(1),P(2)\}}$ for every $x \in \mathcal{A}_{\underline{\{P(1), P(2)\}}}$.}

\begin{center}

\begin{tikzpicture}
 
\node[rectangle,
    draw = black,
    text = black,
    fill = white,
    minimum width = 4.5cm, 
    minimum height = 0.85cm] (1) at (8.5,1.5) {$\mathcal{A}_{\{P(1),P(2)\}}$};

\node[rectangle,
    draw = black,
    text = black,
    fill = white,
    minimum width = 4.5cm, 
    minimum height = 0.85cm] (2) at (4,1.5) {$\mathcal{A}_{\{P(1)\}\underline{\{P(2)\}}}$};

\node[rectangle,
    draw = black,
    text = black,
    fill = white,
    minimum width = 4.5cm, 
    minimum height = 0.85cm] (3) at (4,0.15) {$\mathcal{A}_{\{P(2)\}\underline{\{P(1)\}}}$};

\node[text = black,
    minimum width = 4.5cm, 
    minimum height = 0.67cm] (4) at (2.78,-1.40) {$\mathcal{A}_{\underline{\{P(1),P(2)\}}}$};

\draw (4.85,-0.75) node[anchor=north]{}
  -- (1.775,-0.75) node[anchor=north]{}
  -- (1.775,-2.25) node[anchor=south]{}
  -- (3.1,-2.25) node[anchor=south]{}
  -- cycle; 

\node (5) at (3.85,-1.5) {};
    
\node (6) at (8.63,-1.5) {};

\node (7) at (8.5,-1.63) {};

\node (8) at (8.5,1.2) {};

\draw (5)--(6);

\draw[decoration={markings,mark=at position 1 with
    {\arrow[scale=2,>=To]{>}}},postaction={decorate}] (7)--(8);

\draw [decorate,decoration={brace,amplitude=5pt},xshift=-4pt,yshift=0pt]
(1.5,1) -- (1.5,2) node [black,midway,xshift=-0.6cm] 
{$\mathcal{P}^{(1)}$};

\draw [decorate,decoration={brace,amplitude=5pt},xshift=-4pt,yshift=0pt]
(1.5,-0.35) -- (1.5,0.65) node [black,midway,xshift=-0.6cm] 
{$\mathcal{P}^{(2)}$};

\draw [decorate,decoration={brace,amplitude=5pt},xshift=-4pt,yshift=0pt]
(1.5,-2.30) -- (1.5,-0.66) node [black,midway,xshift=-0.6cm] 
{};

\node[text = black] (9) at (0.05,-1.60) {$\displaystyle \bigcup_{3 \leq i \leq |\mathcal{P}|} \mathcal{P}^{(i)}$};

\end{tikzpicture}

\end{center} 

\end{figure}

\begin{figure}[H]

\caption*{\textbf{Figure 3.2:} A technique for Conjecture 3.3 attempts to find, for any union-closed family $\mathcal{A}$ with $\emptyset \not \in \mathcal{A}$ and $|\mathcal{P}|=2$, a unique $y \in \mathcal{A}_{\{P(1),P(2)\}}$ for every $x \in \mathcal{A}_{\{P(2)\}\underline{\{P(1)\}}}$.}

\begin{center}

\begin{tikzpicture}
 
\node[rectangle,
    draw = black,
    text = black,
    fill = white,
    minimum width = 4.5cm, 
    minimum height = 0.85cm] (1) at (8.5,1.5) {$\mathcal{A}_{\{P(1),P(2)\}}$};

\node[rectangle,
    draw = black,
    text = black,
    fill = white,
    minimum width = 4.5cm, 
    minimum height = 0.85cm] (2) at (4,1.5) {$\mathcal{A}_{\{P(1)\}\underline{\{P(2)\}}}$};

\node[rectangle,
    draw = black,
    text = black,
    fill = white,
    minimum width = 4.5cm, 
    minimum height = 0.85cm] (3) at (4,0.15) {$\mathcal{A}_{\{P(2)\}\underline{\{P(1)\}}}$};

\node (5) at (6.12,0.15) {};  

\node (6) at (8.63,0.15) {};

\node (7) at (8.5,0.03) {};

\node (8) at (8.5,1.2) {};

\draw (5)--(6);

\draw[decoration={markings,mark=at position 1 with
    {\arrow[scale=2,>=To]{>}}},postaction={decorate}] (7)--(8);

\draw [decorate,decoration={brace,amplitude=5pt},xshift=-4pt,yshift=0pt]
(1.5,1) -- (1.5,2) node [black,midway,xshift=-0.6cm] 
{$\mathcal{P}^{(1)}$};

\draw [decorate,decoration={brace,amplitude=5pt},xshift=-4pt,yshift=0pt]
(1.5,-0.35) -- (1.5,0.65) node [black,midway,xshift=-0.6cm] 
{};

\node (9) at (-0.05,0.11) {$\mathcal{P}^{(|\mathcal{P}|)}=\mathcal{P}^{(2)}$};

\end{tikzpicture}

\end{center} 
      
\end{figure}

In both of these techniques, we determine member sets of $\mathcal{A}$ that necessarily belong to $\mathcal{A}_{\{P(1),P(2)\}}$. To do so, we can utilize the facts that $x_1 \in \mathcal{A}_{\{P(1)\}\underline{\{P(2)\}}} \ \land \ x_2 \in \mathcal{A}_{\{P(2)\}\underline{\{P(1)\}}} \implies x_1 \cup x_2 \in \mathcal{A}_{\{P(1),P(2)\}}$ and $y_1 \in \mathcal{A}_{\{P(1),P(2)\}} \ \land \ y_2 \in \mathcal{A} \setminus \mathcal{A}_{\{P(1), P(2)\}} \implies y_1 \cup y_2 \in \mathcal{A}_{\{P(1),P(2)\}}$.

\section*{4. A second strengthening of Conjecture 1.1}

We now consider a strengthening of Conjecture 1.1 that proposes an alternate definition of a finite power set.

\medskip

A family of sets $\mathcal{F}$ is called \textit{separating} if $(x,y \in \bigcup_{F \in \mathcal{F}} F) \ \land \ (\mathcal{F}_{\{x\}} = \mathcal{F}_{\{y\}}) \implies x=y$. In other words, $\mathcal{F}$ is separating if no two distinct elements of its base set are in exactly the same member sets.

\medskip

\noindent \textbf{Conjecture 4.1} (Poonen [6])\textbf{.} \textit{If} $\mathcal{A}$ \textit{is a separating union-closed family and is not a power set, then} $\max_{x \in [n]}\{|\mathcal{A}_{\{x\}}|\}>\frac{|\mathcal{A}|}{2}$\textit{.}

\medskip

Conjecture 4.1 implies Conjecture 1.1 because any finite power set with at least one nonempty member set satisfies Conjecture 1.1, and Conjecture 4.1 states that no other separating union-closed family $\mathcal{A}$ has $\max_{x \in [n]} \{|\mathcal{A}_{\{x\}}|\}$ less than $\frac{|\mathcal{A}|}{2}$. This need not be explicitly stated, however, as the following conjecture also implies Conjecture 1.1.

\medskip

\noindent \textbf{Conjecture 4.2.} \textit{A family of sets} $\mathcal{A}$ \textit{is a finite power set if and only if either} $\mathcal{A}=\{\emptyset\}$ \textit{or} $\mathcal{A}$ \textit{is union-closed and separating with} $\max_{x \in [n]} \{|\mathcal{A}_{\{x\}}|\}=\frac{|\mathcal{A}|}{2}$\textit{.}

\medskip

For families $\mathcal{X}$ and $\mathcal{Y}$, let $\mathcal{X} \bigsqcup \mathcal{Y}=\{X \cup Y \ | \ X \in \mathcal{X} \land Y \in \mathcal{Y} \}$. The following lemma will be used in showing that Conjecture 4.2 implies Conjecture 1.1.

\medskip

\noindent \textbf{Lemma 4.3.} \textit{If} $\mathcal{X}$ \textit{and} $\mathcal{Y}$ \textit{are union-closed families, then} $\mathcal{X} \bigsqcup \mathcal{Y}$ \textit{is also a union-closed family.}

\medskip

\noindent \textit{Proof.} If $\mathcal{X}$ and $\mathcal{Y}$ are union-closed families, then $\mathcal{X} \bigsqcup \mathcal{Y}$ has member sets that are all finite and distinct, and has at least one nonempty member set. To verify the union-closed property for $\mathcal{X} \bigsqcup \mathcal{Y}$, we observe that for all $A$ and $B$ in $\mathcal{X} \bigsqcup \mathcal{Y}$, there are $X$, $X'$ in $\mathcal{X}$ and $Y$, $Y'$ in $\mathcal{Y}$ such that $A = X \cup Y$ and $B = X' \cup Y'$. We must show that $A \cup B = (X \cup Y) \cup (X' \cup Y')$ is in $\mathcal{X} \bigsqcup \mathcal{Y}$. Because $\cup$ is commutative and associative, $(X \cup Y) \cup (X' \cup Y') = (X \cup X') \cup (Y \cup Y')$. $X \cup X' \in \mathcal{X}$ and $Y \cup Y' \in \mathcal{Y}$ because $\mathcal{X}$ and $\mathcal{Y}$ are union-closed. Then $(X \cup X') \cup (Y \cup Y')$ is by definition in $\mathcal{X} \bigsqcup \mathcal{Y}$. Thus, $\mathcal{X} \bigsqcup \mathcal{Y}$ satisfies the union-closed property, and the proof of Lemma 4.3 is complete.

\medskip

\noindent \textbf{Theorem 4.4.} \textit{Conjecture 4.2} $\implies$ \textit{Conjecture 1.1.}

\medskip

\noindent \textit{Proof.} We first note that if a counterexample to Conjecture 1.1 exists, then a separating counterexample to Conjecture 1.1 exists; if a counterexample is not separating, then we can make a separating family on a smaller base set by representing as a single element every grouping of two or more elements (from the original base set) that are contained within exactly the same member sets of the original family. This separating family is also a counterexample as it has the same number of member sets and the same set of element frequencies as of the original family.

\medskip

\noindent Now assume that Conjecture 1.1 is false, and $\tilde{\mathcal{A}}$ is a separating counterexample on a base set $[\tilde{n}]$. We form a union-closed family ${{\mathcal{\tilde{A}}}}'= \tilde{\mathcal{A}} \bigsqcup \{\{z\}, \emptyset\}$ such that $z \not \in [\tilde{n}]$. By Lemma 4.3, ${\tilde{\mathcal{A}}}'$ is union-closed because $\tilde{\mathcal{A}}$ and $\{\{z\}, \emptyset\}$ are both union-closed.

\medskip

\noindent Next, we prove that $\max_{x \in \bigcup_{A \in {\tilde{\mathcal{A}}}'} A}\{|{{\tilde{\mathcal{A}}}'}_{\{x\}}|\}=\frac{|{\tilde{\mathcal{A}}}'|}{2}$.

\medskip

\noindent \textit{Proof.} $({\tilde{\mathcal{A}}}'=\tilde{\mathcal{A}} \cup \{y \cup \{z\} \ | \ y \in \tilde{\mathcal{A}} \}) \ \land \ (z \not \in [\tilde{n}]) \implies |{\tilde{\mathcal{A}}}'|=2|\tilde{\mathcal{A}}|$. Also, $\forall x \in [\tilde{n}] \ (({{\tilde{\mathcal{A}}}'}_{\{x\}}=\tilde{\mathcal{A}}_{\{x\}} \cup \{y \cup \{z\} \ | \ y \in \tilde{\mathcal{A}}_{\{x\}} \}) \ \land \ (z \not \in [\tilde{n}]) \implies |{{\tilde{\mathcal{A}}}'}_{\{x\}}|=2|\tilde{\mathcal{A}}_{\{x\}}|)$. Therefore, $\forall x \in [\tilde{n}] \ \frac{|{{\tilde{\mathcal{A}}}'}_{\{x\}}|}{|{\tilde{\mathcal{A}}}'|}=\frac{|\tilde{\mathcal{A}}_{\{x\}}|}{|\tilde{\mathcal{A}}|} < \frac{1}{2}$, making all elements, except for $z$, in $\bigcup_{A \in {\tilde{\mathcal{A}}}'} A=[\tilde{n}]\cup \{z\}$ have frequency in ${\tilde{\mathcal{A}}}'$ less than $\frac{|{\tilde{\mathcal{A}}}'|}{2}$. $|{{\tilde{\mathcal{A}}}'}_{\{z\}}|=\frac{|{\tilde{\mathcal{A}}}'|}{2}$, as $|{{\tilde{\mathcal{A}}}'}_{\{z\}}|+|{{\tilde{\mathcal{A}}}'}_{\underline{{\{z\}}}}|=|{{\tilde{\mathcal{A}}}'}|$ and $({\tilde{\mathcal{A}}}'=\tilde{\mathcal{A}} \cup \{y \cup \{z\} \ | \ y \in \tilde{\mathcal{A}} \}) \ \land \ (z \not \in [\tilde{n}]) \implies |{{\tilde{\mathcal{A}}}'}_{\{z\}}|=|{{\tilde{\mathcal{A}}}'}_{\underline{{\{z\}}}}|$. Thus, $z$ has the highest frequency in $\mathcal{{\tilde{\mathcal{A}}}'}$ among elements of $\bigcup_{A \in {\tilde{\mathcal{A}}}'} A$, and $\max_{x \in \bigcup_{A \in {\tilde{\mathcal{A}}}'} A}\{|{{\tilde{\mathcal{A}}}'}_{\{x\}}|\}=\frac{|{\tilde{\mathcal{A}}}'|}{2}$.

\medskip

\noindent We also have that ${\tilde{\mathcal{A}}}'$ is separating. To prove this, we observe that $\{x,y\} \in \binom{[\tilde{n}]}{2} \implies \tilde{\mathcal{A}}_{\{x\}} \neq \tilde{\mathcal{A}}_{\{y\}}$, as $\tilde{\mathcal{A}}$ itself is separating. It follows that $\{x,y\} \in \binom{[\tilde{n}]}{2} \implies {\tilde{\mathcal{A}}_{\{x\}}'} \neq {\tilde{\mathcal{A}}'}_{\{y\}}$. Now, $\forall x \in [\tilde{n}] \ (\forall A \in \tilde{\mathcal{A}}_{\{x\}} \ (A \in {\tilde{\mathcal{A}}}' \ \land \ z \not \in A)) \implies {\tilde{\mathcal{A}}'}_{\{x\}} \neq {\tilde{\mathcal{A}}'}_{\{z\}}$. Thus, no two distinct elements of the base set of ${\tilde{\mathcal{A}}}'$ are in exactly the same member sets of ${\tilde{\mathcal{A}}}'$.

\medskip

\noindent ${\tilde{\mathcal{A}}}'$ is not a finite power set, as all elements of $[\tilde{n}]$ have nonzero frequency in ${\tilde{\mathcal{A}}}'$ less than $\frac{|{\tilde{\mathcal{A}}}'|}{2}$.

\medskip

\noindent Hence, ${\tilde{\mathcal{A}}}'$ is a counterexample to Conjecture 4.2, as ${\tilde{\mathcal{A}}}' \neq \{ \emptyset \}$, ${\tilde{\mathcal{A}}}'$ is union-closed and  separating with $\max_{x \in \bigcup_{A \in {\tilde{\mathcal{A}}}'} A}\{|{{\tilde{\mathcal{A}}}'}_{\{x\}}|\}=\frac{|{\tilde{\mathcal{A}}}'|}{2}$, and ${\tilde{\mathcal{A}}}'$ is not a finite power set. We have that $\neg$ (Conjecture 1.1) $\implies \neg$ (Conjecture 4.2). Therefore, Conjecture 4.2 $\implies$ Conjecture 1.1.

\medskip

A motivation for decoupling Conjecture 4.2 from Conjecture 4.1 is to emphasize possible equivalence with Conjecture 1.1. Conjecture 4.2 implies Conjecture 1.1, but it remains to be shown whether Conjecture 1.1 also implies Conjecture 4.2.

\end{document}